\documentclass{amsart}

\usepackage{amsopn, amsthm, amsgen, amscd,amsmath, amssymb}
\usepackage[all]{xy}
\usepackage{hyperref}

\newcommand{\N}{{\mathbb N}}
\newcommand{\Z}{{\mathbb Z}}
\newcommand{\R}{{\mathbb R}}
\newcommand{\C}{{\mathbb C}}
\newcommand{\Q}{{\mathbb Q}}
\newcommand{\tr}{\mathop{\rm tr}}
\newcommand{\dem}{{\em Proof: \;}}
\newcommand{\fdem}{\hfill $\square$}

\newtheorem{theorem}{Theorem}

\theoremstyle{plain}
\newtheorem{teo}{Theorem}[section]
\newtheorem{lema}[teo]{Lemma}
\newtheorem{cor}[teo]{Corollary}
\newtheorem{prop}[teo]{Proposition}

\theoremstyle{definition}
\newtheorem{defi}{Definition}[section]

\theoremstyle{remark}
\newtheorem{rmrk}{Remark}[section]

\title[]{Teichm\"uller theory of the universal hyperbolic lamination}

\author{J. M. Burgos}
\address{Departamento de Matem\'aticas, Centro de Investigaci\'on y de Estudios Avanzados, Av. Instituto Polit\'ecnico Nacional 2508, Col. San Pedro Zacatenco, C.P. 07360 Ciudad de M\'exico, M\'exico}
\email{burgos@math.cinvestav.mx}
\author{A. Verjovsky}
\address{Instituto de Matem\'{a}ticas -- Unidad Cuernavaca, Universidad
Nacional Aut\'{o}noma de M\'{e}xico, Av. Universidad S/N, C.P. 62210
Cuernavaca, Morelos, M\'{e}xico}
\email{alberto@matcuer.unam.mx}



\begin{document}

\title[]{Teichm\"uller theory of the universal hyperbolic lamination}

\begin{abstract}
We construct an Ahlfors-Bers complex analytic model for the Teichm\"uller space of the universal hyperbolic lamination (also known as Sullivan's Teichm\"uller space) and the renormalized Weil-Petersson metric on it as an extension of the usual one. In this setting, we prove that Sullivan's Teichm\"uller space is K\"ahler isometric biholomophic to the space of continuous functions from the profinite completion of the fundamental group of a compact Riemann surface of genus greater than or equal to two to the Teichm\"uller space of this surface; i.e. We find natural K\"ahler coordinates for the Sullivan's Teichm\"uller space. This is the main result. As a corollary we show the expected fact that the Nag-Verjovsky embedding is transversal to the Sullivan's Teichm\"uller space contained in the universal one.
\end{abstract}

\vspace{2cm}

\maketitle

\section{Introduction}

In his celebrated paper, which connects the universalities of Milnor-Thurston, Feigenbaum with the Ahlfors-Bers theory  \cite{Universalities}, Dennis Sullivan introduces in Example 4, Sect 1 a lamination (or solenoidal surface) by Riemann surfaces as the inverse limit of the tower of finite covers of a surface $\Sigma_2$ of genus two. This is
in the spirit of Grothendieck's algebraic universal covering, and it is called the \emph{Universal Hyperbolic Lamination}. This lamination $\Sigma_{\infty}$ is
a principal fiber bundle $p:\Sigma_{\infty}\to\Sigma_2$  with fiber the group $\hat{G}$ which is profinite completion of the fundamental group $G$ of $\Sigma_2$. This is a topological nonabelian group which is homeomorphic to the Cantor set. Since every surface of genus greater or equal to 2 covers $\Sigma_2$, the group $\hat{G}$ is also the profinite completion of the fundamental group of any  surface of genus $g\geq2$.

In the same paper, Sect 2, Sullivan defined the Teichm\"uller space of this 2-dimensional lamination as follows:

\begin{quote}
``A \textit{complex structure} on a lamination $L$ is a maximal covering of $L$ by lamination charts (disk)$\times$(transversal) so that overlap homeomorphisms are complex analytic in the disk direction. Two complex structures are \textit{Teichm\"uller equivalent} if they are related by a homeomorphism which is homotopic to the identity through leaf preserving continuous mappings of $L$. The set of classes is called the \textit{Teichm\"uller set} $T(L)$."
\end{quote}
Given a complex structure $J$ on $\Sigma_2$, using the projection $p:\Sigma_{\infty}\to\Sigma_2$ we can pull-back the complex structure to the leaves of the lamination to give a complex structure $p^*(J)$ to $\Sigma_\infty$. One has a canonical inclusion of the Teichm\"uller of $\Sigma_2$ into the Teichm\"uller of $\Sigma_\infty$.

In section \ref{Teichmuller_Space} we give an equivalent definition to the Teichm\"uller space, and in section \ref{Limit_automorphic_differentials} we give an Ahlfors-Bers complex analytic model for this space based on the notion limit-automorphic Beltrami differentials defined in the same section. In section \ref{Ren_Weil_Petersson} we show that, in contrast to the universal Teichm\"uller space, the Sullivan's version admits a well defined Weil-Petersson metric extension: the Renormalized Weil-Petersson metric. As in the classical finite dimensional case, this is a K\"ahler metric. However, this metric is not induced by the canonical embedding $T(\Sigma_{\infty})\hookrightarrow T(1)$. Nevertheless, the transversality\footnote{We don't assume that the sum of the respective tangent spaces at the intersection is the whole space. We just say that the intersection is transversal if the sum of the tangent spaces at the intersection point is direct. This can be done because we are not interested in any intersection invariant in this work.} respect to the Nag-Verjovsky embedding \cite{NV} still holds:
$$Diff^{+}(S^{1})/PSL(2,\R)\pitchfork T(\Sigma_\infty)$$
This is proved in Theorem \ref{Trans_Nag_Verj}. This result was expected considering the identity in Proposition \ref{Limit_Teichm}:
$$T\left(\Sigma_{\infty}\right)= \overline{\lim_{\substack{\longrightarrow \\ G'<G \\ [G':G]<\infty}} T\left(\Sigma_{G'}\right)}
\quad \rm{(closure\, of\, the\, direct\, limit)}
$$

Now, with these tools at hand, we prove the Theorem A:

\begin{theorem}\label{main_Th_Intro}
Consider a genus $g\geq 2$ compact Riemann surface $\Sigma_{g}$ and its universal hyperbolic lamination $\Sigma_{\infty}$ fibering over the surface. Respect to the renormalized Weil-Petersson metric, there is a K\"ahler isometric biholomorphism between the Sullivan's Teichm\"uller space and the space of continuous functions from the pro-finite completion  $ \widehat{\pi_{1}(\Sigma_{g})}$ of the fundamental group to the Teichm\"uller space of the surface:
$$C\left( \widehat{\pi_{1}(\Sigma_{g})}, T(\Sigma_{g})\right) \cong T(\Sigma_{\infty})$$
where $\Sigma_{\infty}$ is the universal hyperbolic lamination and $T(\Sigma_{\infty})$ is Sullivan's Teichm\"uller space. Moreover, the biholomorphism is functorial respect to the inverse limit maps of the pro-finite completion. The K\"ahler structure in the space $C\left( \widehat{\pi_{1}(\Sigma_{g})}, T(\Sigma_{g})\right)$  is induced by the K\"ahler structure given by the Weil-Petersson metric in $T(\Sigma_{\infty})$.
\end{theorem}

We can think of this result as finite-dimensional valued K\"ahler coordinates on Sullivan's Teichm\"uller space. Because of the functoriality respect to the inverse limit maps, we have Theorem B:

\begin{theorem}
There is a discrete fiber holomorphic K\"ahler isometric branched covering (in the orbifold sense):
$$\mathcal{M}_{g}^{n}\twoheadrightarrow \mathcal{M}_{n(g-1)+1}$$
Moreover, this covering factors through the alternate product $Alt^{n}(\mathcal{M}_{g})$.
\end{theorem}

Note that this result is not in contradiction with the rigidity result in \cite{Aramayona_Souto} \footnote{
Consider a holomorphic map between moduli spaces:
$$\phi:\mathcal{M}_{g, k}\rightarrow \mathcal{M}_{g', k'}$$
If $g'\leq  2g-2$ and $g\geq 6$, then $\phi$ is either constant or a forgetful map. In the later case, $g'=g$ and $k'\leq k$.}. The simplest case is\footnote{Is there any explicit way of describing a genus three Riemann surface by two genus two surfaces? Do not think about nodal curves, we are not in the Deligne-Mumford boundary.}:
$$\mathcal{M}_{2}\times\mathcal{M}_{2}\twoheadrightarrow \mathcal{M}_{3}$$

In section \ref{Siegel_Functions} we define the space of Siegel functions as a substitute of the Segal disk for the Sullivan's Teichm\"uller space.

\section{Complex analytic Teichm\"uller theory in a nutshell}

Consider a Riemann surface $\Sigma$. How we deform its complex structure? Consider the Poincar\'e-Koebe uniformization $\Delta\rightarrow \Sigma$ and the representation of $G:= \pi_{1}(\Sigma)$ as a Fuchsian group:
$$\alpha:G\rightarrow Isom^{+}(\Delta)$$

\begin{defi}
\ \\
\begin{itemize}
\item An $L_{\infty}$ section of the line bundle $\mu\in \overline{\omega}\otimes\omega^{*}$ where $\omega$ is the canonical line bundle of $\Delta$ will be just called a \textit{differential}; i.e. It is a $(-1, 1)$ form and its local expression is  $\mu= h\otimes d\bar{z}\otimes \partial_{z}$ where $h\in L_{\infty}(\Delta)$.
\item A differential $\mu$ is a $G$\textit{-automorphic} differential if it is a differential and:
$$\alpha(g)^{*}\mu= \mu\ \ \ \forall\ g\in G$$
The space of $G$-automorphic differentials will be denoted by $L_{\infty}(G)$.
\item A differential $\mu$ is a \textit{Beltrami differential} if it is a differential and $||\mu||_{\infty}<1$. These differentials are the deformation parameters. 
\end{itemize}
\end{defi}

\begin{rmrk}
The pullback of a differential in $\Sigma$ by the uniformization map is a $G$-automorphic differential in the Poincar\'e disk.
\end{rmrk}

How a Beltrami differential, the deformation parameter, actually realizes a deformation? A Beltrami differential can be seen as an $\infty$-measurable field of homothetic classes of ellipses: The relation between the major and minor axis is $K=(1+ |\mu|)/(1- |\mu|)$ and the major axis is rotated by $arg(\mu)/2$.

Consider a Beltrami differential on the complex plane such that $\left(\bar{z}^{-1}\right)^{*}\mu=\mu$; i.e. It is symmetric respect to the unit circle. Consider the Ahlfors-Bers equation:
$$\partial_{\bar{z}}f= \mu\ \partial_{z}f$$
Is there any solution to this equation on the disk $\Delta$? Geometrically, Is there a map $f$ on the disk straightening all the infinitesimal ellipses into infinitesimal circles? 

The following is the Ahlfors-Bers existence Theorem adapted to our situation \cite{Imayoshi}:

\begin{teo}
Consider a Beltrami differential on the complex plane such that $\left(\bar{z}^{-1}\right)^{*}\mu=\mu$. There are quasiconformal homeomorphisms solutions to the Ahlfors-Bers equation $\partial_{\bar{z}}f= \mu\ \partial_{z}f$. Moreover, these solutions uniquely extends to a quasisymmetric homeomorphism on the boundary and there is a unique solution $f^{\mu}$ fixing $1,\ i$ and $-1$.
\end{teo}

If $\mu=0$, then $\partial_{\bar{z}}f^{\mu}=0$ and by the Weil Lemma, $f^{\mu}$ is holomorphic. Define the pair $(\Sigma_{\mu}, [f^{\mu}])$ by the following commutative diagram:
$$\xymatrix{\Delta \ar[rr]^{f^{\mu}} \ar[d]& & \Delta \ar[d] \\ \Sigma\cong\Delta/G \ar[rr]^{[f^{\mu}]} & & \Sigma_{\mu}\cong\Delta/G_{\mu} }$$
where $G_{\mu}$ is the group of conjugated elements $f\gamma f^{-1}$ such that $\gamma\in G$. The group $G_{\mu}$ is also Fuchsian hence $\Sigma_{\mu}$ is a Riemann surface. The surface $\Sigma_{\mu}$ will be called the $\mu$-deformed Riemann surface.

Are we really deforming? Is there any redundancy in the parameters? We say that the pairs $(\Sigma_{\mu}, [f^{\mu}])$ and $(\Sigma_{\nu}, [f^{\nu}])$ are equivalent if there is a biholomorphism $h$ such that the following diagram commutes up to homotopy:

$$\xymatrix{ & \Sigma\ar[dl]_{[f^{\mu}]} \ar[dr]^{[f^{\nu}]} & \\
\Sigma_{\mu}\ar[rr]^{h} & & \Sigma_{\nu}}$$

The equivalence classes constitutes the \textit{Teichm\"uller space} $T(\Sigma)$ of compact Riemann surfaces of genus $g\geq 2$.
This is a fine moduli space and the universal family will be called the \textit{Teichm\"uller universal curve}. The coarse moduli space $\mathcal{M}(\Sigma)$ of compact Riemann surfaces is the quotient by the \textit{Mapping Class Group}:

$$\mathcal{M}(\Sigma)= T(\Sigma)/MCG(\Sigma)$$
$$MCG(\Sigma):=Homeo^{+}(\Sigma)/Homeo_{0}(\Sigma)\cong Aut\left(\pi_1(\Sigma)\right)$$

\begin{prop}
The pairs $(\Sigma_{\mu}, [f^{\mu}])$ and $(\Sigma_{\nu}, [f^{\nu}])$ are equivalent if and only if the quasisymmetric homeomorphisms on the boundary $f^{\mu}|_{\partial\Delta}$ and $f^{\nu}|_{\partial\Delta}$ are equal.
\end{prop}

The above condition defines an equivalence relation $\sim$ in the space of $G$-automorphic Beltrami differentials $L_{\infty}(G)_{1}$. We have the following Ahlfors-Bers complex analytic model for the Teichm\"uller space:
$$T(\Sigma)\cong L_{\infty}(G)_{1}/\sim$$
The tangent space at a point $[\mu]$ is $L_{\infty}(G)$ quotiented by the induced equivalence relation. Define the complex analytic structure as follows: $J[\nu]_{\mu}:=[i\nu]_{\mu}$ where $[\nu]_{\mu}$ is a tangent vector at a point $[\mu]$. The following deep Theorem characterize the complex analytic structure of the Teichm\"uller space:

\begin{teo}
$T(\Sigma_{g})$ is a complex domain of complex dimension $3g-3$.
\end{teo}

The unit ball of Beltrami differentials $L_{\infty}(\Delta)_{1}$ has a group structure such that $f^{\mu*\nu^{-1}}= f^{\mu}\circ (f^{\nu})^{-1}$ where the inverse $\nu^{-1}$ is respect to this product. Explicitly:
$$\mu*\nu^{-1}:=\frac{\mu-\nu}{1-\mu\bar{\nu}}\frac{(f^{\nu})_{z}}{(\overline{f^{\nu}})_{\bar{z}}}\circ (f^{\nu})^{-1}$$
This product projects to $T(1)$: $[\alpha]*[\beta]:=[\alpha*\beta]$. Respect to this product, $T(1)$ is not a topological group: Right translations $R_{[\mu]}:T(1)\rightarrow T(1)$ such that $R_{[\mu]}(\nu):= \nu*\mu$ act by biholomorphic automorphisms of $T(1)$ while the left translations are not even continuous (\cite{Leh}, Sect.III.3.4).

The Weil-Petersson metric on $T(\Sigma)$ is defined at the origin $(\Sigma, id)$ by the following expression on representatives:
$$WP_{G}(\mu, \nu):= \int_{F}dxdy\int_{\Delta}dx'dy'\ \frac{\mu(z)\overline{\nu}(\zeta)}{(1-z\bar{\zeta})^{4}}$$
such that $F$ is a fundamental domain of the action by $G$ and $\mu$ and $\nu$ are tangent vectors at the origin; i.e. elements of $L_{\infty}(G)$. Translating this metric by the right translations $R_{[\mu]}$ we have a K\"ahler metric on the Teichm\"uller space respect to the complex analytic structure defined before \cite{Takhtajan}.

By the Ahlfors-Beurling extension Theorem, every quasisymmetric homeomorphism on the circle extends to a quasiconformal homeomorphism on the disk. Then:
$$T(1):= L_{\infty}(\Delta)_{1}/\sim\cong QS(S^{1})/PSL(2,\R)$$
where $T(1)$ is the \textit{universal Teichm\"uller space}. It is universal in the sense that it contains all the finite dimensional Teichm\"uller spaces:
$$T(\Sigma_g)\subset T(1)$$

Because every diffeomorphism of the circle is quasisymmetric, we have the Nag-Verjovsky embedding of the coadjoint orbit $Diff^{+}(S^{1})/PSL(2,\R)$ into the universal Teichm\"uller space:
$$Diff^{+}(S^{1})/PSL(2,\R)\hookrightarrow T(1)$$
In \cite{NV} it is proved that this map is a K\"ahler isometric complex analytic embedding respect to the natural K\"ahler structure of the coadjoint orbit and the K\"ahler structure of $T(1)$ with the Weil-Petersson metric\footnote{The Weil-Petersson metric on $T(1)$ is defined only for tangent vectors whose associated vector fields on the circle are $C^{3/2+\varepsilon}$}. It defines a holomorphic foliation of $T(1)$ by $Diff^{+}(S^{1})/PSL(2,\R)$ leaves. This is another way to characterize the complex analytic structure of $T(1)$.

The Nag-Verjovsky embedding is transversal to every finite dimensional Teichm\"uller space:
$$Diff^{+}(S^{1})/PSL(2,\R)\pitchfork T(\Sigma_g)$$
This is an infinitesimal form of Mostow rigidity.

\section{Laminated uniformization}

Consider a Riemann surface $\Sigma$\footnote{The following construction works over any cocompact finitely generated Fuchsian group as well.}. By Poincar\'e-Koebe uniformization we have a faithful representation of the fundamental group $\pi_{1}(\Sigma)$ as a Fuchsian group\footnote{Because $\Sigma$ is a compact surface, the Fuchsian group consists of hyperbolic transformations only.}:
$$\alpha: \pi_{1}(\Sigma)\rightarrow Isom^{+}(\Delta)\cong\text{PSL}(2,\R)$$
where $\Delta$ is the hyperbolic disk and $Isom^{+}(\Delta)$ is the group of positively oriented hyperbolic isometries. We will denote by $G$ the image of the fundamental group by $\alpha$; i.e. $G:= \alpha(\pi_{1}(\Sigma))$. We recover the Riemann surface $\Sigma$ as the quotient of the disk $\Delta$ by the $G$-action:
$$\Sigma\simeq\Delta/G$$

Consider a finite unramified holomorphic connected covering $\Sigma'$ of the Riemann surface $\Sigma$ with finite fiber $\mathcal{F}_{n}=\{1,2\ldots n\}$:

\begin{equation}\label{covering}
\xymatrix{
\mathcal{F}_{n}  \ar@{^{(}->}[r] & \Sigma' \ar[d]_{p}^{n:1} \\
& \Sigma}
\end{equation}

Denote by $\rho$ the monodromy representation of the fundamental group $\pi_{1}(\Sigma)$ respect to the covering $p:\Sigma'\rightarrow\Sigma$ \footnote{This is well defined because the covering is unramified. The action on the fiber $\mathcal{F}_{n}$ is a right action.}:
\begin{equation}\label{monodromy}
\rho:\pi_{1}(\Sigma)\rightarrow Sym_{n}^{op}
\end{equation}
where $Sym_{n}^{op}$ is the opposite group of the symmetric group of permutations of the set $\mathcal{F}_{n}$. Because the covering $\Sigma'$ is connected, the representation $\rho$ is transitive.

Conversely, given a transitive representation \eqref{monodromy} we recover the covering \eqref{covering} as follows: Consider the tensor product representation of the representations $\rho^{-1}$ and $\alpha$ as before \footnote{We take the monodromy inverse because we look for the following equivalence relation (Recall that $\rho$ is a right action):
$$((x)\rho(g),z)\sim (x,\alpha(g)(z))$$}:
$$\rho^{-1}\otimes\alpha:\pi_{1}(\Sigma)\rightarrow Isom^{+}(\mathcal{F}_{n}\times\Delta)$$
The trivial bundle $p_{2}:\mathcal{F}_{n}\times\Delta\rightarrow\Delta$ is equivariant respect to the action by the fundamental group; i.e. For every element $g\in \pi_{1}(\Sigma)$ the following diagram commutes:
\begin{equation}\label{equivariance_finite}
\xymatrix{\mathcal{F}_{n}\times\Delta \ar[d]_{p_{2}} \ar[rrr]^{\rho^{-1}\otimes\alpha(g)} & & & \mathcal{F}_{n}\times\Delta \ar[d]^{p_{2}} \\
\Delta \ar[rrr]^{\alpha(g)} & & & \Delta}
\end{equation}
Since the representation $\alpha$ is proper, free and discontinuous so is the representation $\rho^{-1}\otimes\alpha$ hence the quotient of the trivial bundle by the equivariant action is a finite unramified holomorphic covering \eqref{covering} such that\footnote{Because the Fuchsian groups are isomorphic, we make an abuse of notation and denote them by the same letter $G$.}:
$$\Sigma'=(\mathcal{F}_{n}\times\Delta)/G$$
is a Riemann surface and $p$ is defined such that $p\left(G\cdot(m,z)\right)= p_{2}(m,z)=z$. The monodromy of the covering is the original representation $\rho$ and because it is transitive the covering is connected.

Because the above constructions were functorial, we have proved the following Theorem:

\begin{teo}\label{Category_Equivalence}
Consider a Riemann surface $\Sigma$. Then, the category of finite (unramified) holomorphic connected coverings of the Riemann surface $\Sigma$ is equivalent to the category of finite irreducible representations $\mathcal{F}_{n}=\{1,2\ldots n\}$ of the fundamental group $\pi_{1}(\Sigma)$.
\end{teo}

Consider the covering $p:\Sigma'\rightarrow\Sigma$ as before and the Fuchsian group $G'$ of the Riemann surface $\Sigma'$; i.e. a discrete subgroup of $Isom^{+}(\Delta)$ such that:
$$\Sigma'=\Delta/G'$$
The group $G'$ is a finite index subgroup of $G$ and as a set, the fiber $\mathcal{F}_{n}$ is in natural one to one correspondence with the coset set $G'\backslash G$. In particular the index of the subgroup equals the number of sheets of the covering. If the subgroup $G'$ is normal, then the fiber is actually a group and the monodromy \eqref{monodromy} is the canonical $G$ group left action on $G'\backslash G$. Conversely, every finite index normal subgroup $G'$ of $G$ defines the monodromy \eqref{monodromy} as the canonical left group action of $G$ in the quotient $G'\backslash G$ where the fiber $\mathcal{F}_{n}$ is identified with this quotient as a set.

\begin{cor}\label{Category_Equivalence_II}
Consider a Riemann surface $\Sigma$. Then, the category of finite (unramified) holomorphic connected regular\footnote{We say a covering is \textit{regular} if given a closed path in the base space all of its preimages are either closed or open paths.} coverings of the Riemann surface $\Sigma$ is equivalent to the category of finite index normal subgroups of the fundamental group $\pi_{1}(\Sigma)$.
\end{cor}

The equivalence functor is the following: Given a finite index subgroup $G'$ of $G$, the corresponding Riemann is surface is:
\begin{equation}\label{Sigma_G_prime}
\Sigma_{G'}:=\left((G'\backslash G)\times \Delta\right)/G
\end{equation}
where the action of $G$ is the diagonal action:
\begin{equation}\label{Action}
g\cdot(h,z)=\left( h\cdot g^{-1}, \alpha(g)(z)\right)
\end{equation}

\begin{table}
\begin{center}
    \begin{tabular}{ | l | p{7cm} |}
    \hline
    Coverings of $\Sigma$ & Subgroups of $\pi_{1}(\Sigma)$ \\ \hline \hline
    Trivial & $\pi_{1}(\Sigma)$  \\ \hline
    Universal & Trivial subgroup $\{e\}$  \\ \hline
    Finite & Finite index  \\ \hline
    Number of sheets & Subgroup index  \\ \hline
    Regular & Normal  \\ \hline
    \end{tabular}
\end{center}
\caption{Dictionary between coverings and subgroups.}\label{Dictionary}
\end{table}

Table \ref{Dictionary} shows a dictionary between coverings and groups and we suggest the classical reference \cite{Seifert} for the proofs. The categories of the previous Theorems are actually inverse systems. The inverse limit of irreducible finite index normal subgroups of $G$ is its profinite completion group $\hat{G}$ \footnote{Because the subsystem of normal finite index subgroups is cofinal in the system of finite index subgroups, this definition coincides with the usual one. This definition underlines the fact that the inverse limit $\hat{G}$ is a group and the inverse limit morphisms $\hat{G}\rightarrow G'\backslash G$ are actually group morphisms.}:
$$\hat{G}=\lim_{\substack{\longleftarrow \\ G'\triangleleft G \\ [G':G]<\infty}} G'\backslash G$$
The collection of finite index subgroups of $G$ is a neighborhood system of the identity and by translation it defines a topology on $G$ whose completion is the group $\hat{G}$ just defined. As a topological space, the group $\hat{G}$ is a compact, perfect and totally disconnected Hausdorff space; i.e. a Cantor set. Because the Fuchsian group $G$ is residually finite, the completion is a group extension and we have a dense immersion:
\begin{equation}\label{dense_inmersion_group}
G\hookrightarrow \hat{G}
\end{equation}

Consider the cofinal inverse system $\left(A_{n}\right)_{n\in\N}$ such that $A_{n}$ is the intersection of all subgroups of index $n$ or less in $G$. Because the group $G$ is finitely generated, there must be a finite amount of subgroups of a given index hence the normal subgroups $A_{n}$ are of finite index as well. Define the following valuation $val:G\rightarrow \N\cup\{\infty\}$ such that:
\begin{equation}\label{valuation}
val(g):=max\{n\in\N\ \ |\ \ g\in A_{n}\}
\end{equation}
if $g$ is not the neutral element $e$ and $val(e):=\infty$. Define the translation invariant metric $d$ on the group $G$ such that:
\begin{equation}\label{distance}
d(g,h):=e^{-val(g^{-1}h)}
\end{equation}
This metric is a non-archimedean metric\footnote{This metric verifies the stronger condition:
$$d(g,h)\leq \max\{d(g,f),d(f,h)\}$$
for every $f,g$ and $h$ in the group $\hat{G}$. This condition implies the usual triangle inequality.} whose induced topology is the one just defined and whose completion is isomorphic to the previous profinite completion $\hat{G}$.

The inverse limit of the category of unramified holomorphic connected coverings of the Riemann surface $\Sigma$ is the \textit{universal hyperbolic lamination} $\Sigma_{\infty}$ \cite{Universalities}. The functoriality of the previous construction and the fact that the regular coverings system is cofinal yields the following identity:
\begin{equation}\label{laminated_unif}
\Sigma_{\infty}=(\hat{G}\times\Delta)/G
\end{equation}
where the action on $\hat{G}\times\Delta$ by the Fuchsian group $G$ is the diagonal action as before.

Identity \eqref{laminated_unif} will be called the \textit{laminated uniformization}\footnote{In the paper \cite{Odden}, this is called a $G$\textit{-tagged solenoid}.} of the universal hyperbolic lamination $\Sigma_{\infty}$ and the space $\hat{G}\times\Delta$ will be called the \textit{laminated Poincar\'e disk}.

Consider the embedding $\iota: \Delta\hookrightarrow \hat{G}\times \Delta$ such that $\iota(z)=(e,z)$ where $e$ is the neutral element of the group $\hat{G}$. The immersion:
$$\xymatrix{ \Delta \ar@{^{(}->}[rr]^{\iota} & & \hat{G}\times \Delta \ar[rr]^{\pi} & & \Sigma_{\infty}=(\hat{G}\times\Delta)/G}$$ 
will be called the \textit{baseleaf} of the lamination $\Sigma_{\infty}$. Because the morphism \eqref{dense_inmersion_group} is a dense immersion so is the baseleaf. It is clear from the identity \eqref{laminated_unif} that the leaf space of the lamination $\Sigma_{\infty}$ is the group $\hat{G}/G$.

From the definition of the baseleaf, we have that the initial coarser topology on the Poincar\'e disk is given by the following basic open set:
$$V(G',U):=\bigcup_{g\in G'}\ g(U)$$
where $U$ is an open set of the Poincar\'e disk and $G'$ is a finite index subgroup of $G$. If the open set $U$ is contained in some fundamental domain of the action by $G$, then the corresponding basic set is a disjoint union. The disk with this coarser topology will be called the \textit{embedded disk} and will be denoted by $\Delta_{Emb}$\footnote{In contrast with the baseleaf immersion, $\Delta_{Emb}\hookrightarrow \Sigma_{\infty}$ is an embedding.}.

\begin{defi}
A continuous function $f:\Delta\rightarrow X$ where $X$ is a topological space will be called \textit{limit-automorphic} if the function $f:\Delta_{Emb}\rightarrow X$ is continuous\footnote{This is the generalization of the usual definition: $f:\R\rightarrow\R$ is \textit{limit-periodic} if for every $\varepsilon>0$ there is a natural number $N$ such that $|f(x+2\pi n)-f(x)|<\varepsilon$ if $n\in N\Z$ for every real number $x$. This is a particular case of the \textit{almost periodic} property \cite{Bohr} where the relatively dense sets are restricted to be the subgroups $N\Z$. There is a beautiful discussion of almost periodic functions in the context of constructive mathematics in \cite{Brom}. An interesting discussion relating limit-periodic functions, solenoids and adding machines can be found in \cite{Bell}}.
\end{defi}

\begin{lema}
The following are equivalent:
\begin{enumerate}
\item $f$ is limit-automorphic.
\item There is a net of automorphic functions $\left(f_{G'}\right)$ ranging over the finite index subgroups of $G$ such that $f_{G'}$ is $G'$-automorphic  and the net uniformly converges to $f$.
\item There is a unique continuous function $f_{\infty}$ on the lamination $\Sigma_{\infty}$ such that $(\pi\iota)^{*}(f_{\infty})= f$.
\end{enumerate}
\end{lema}
\dem The proof is verbatim to the one in Lemma \ref{TFAE}.
\fdem

\section{Beltrami differentials and Teichm\"uller space}

Consider the normalized Haar measure $\eta$ on the group $\hat{G}$. By Caratheodory's measure extension Theorem, there is a finite Borel measure $m$ on the embedded disk $\Delta_{Emb}$ such that:
$$m\left(V(G',U)\right)= \eta(\overline{G'})\lambda(U)$$
where $\lambda$ is the Borel measure induced by the hyperbolic metric of the Poincar\'e disk and the closure is respect to the profinite topology of the group $G$ \footnote{The measure $\lambda$ is invariant under hyperbolic isometries hence the definition of the measure $m$ on basic open sets $V(G',U)$ is well defined.}. The measure $m$ of the disk is the hyperbolic area of a fundamental domain of the action $G$ hence is finite for $G$ is cocompact.

The pushforward of the diagonal action invariant measure $\eta\times m$ on the laminated disk by the laminated uniformization $\pi:= \hat{G}\times\Delta\rightarrow \Sigma_{\infty}$ defines a finite measure $m_{\infty}$ on the lamination $\Sigma_{\infty}$ \footnote{The reader may be tempted to define $m_{\infty}$ as the pushforward of the $G$-invariant measure $\eta\times \lambda$. But this measure would be trivial in the sense that it would take the values $0$ and $\infty$ only. In particular it wouldn't be even $\sigma$-finite.}. If we desintegrate this measure along the fibration over a trivializing open set, it becomes the product measure $\eta\times\lambda$. Again, the measure $m_{\infty}$ of the lamination is the hyperbolic area of a fundamental domain of the action by $G$ hence is finite for $G$ is cocompact.

\subsection{Deformations and Teichm\"uller space}\label{Teichmuller_Space}

A leaf preserving homeomorphism from the laminated Poincar\'e disk into itself will be called canonical quasiconformal if it is canonical quasiconformal leaf by leaf. For every differential $\mu\in L_{\infty}(\Delta)$ define the extension $\tilde{\mu}$ such that:
$$\tilde{\mu}(z):=(1/\bar{z})^{*}(\mu)(z)=\mu\left(\frac{1}{\bar{z}}\right)\frac{\bar{z}^{2}}{z^{2}}$$
for every $z\in \C$ such that $|z|>1$. See that $||\mu||_{\infty}=||\tilde{\mu}||_{\infty}$.

\begin{lema}
Consider a function $\zeta: \hat{G}\rightarrow L_{\infty}(\Delta)_{1}$. The Ahlfors-Bers equation:
\begin{equation}\label{Ahlfors_Bers}
\partial_{\bar{z}}f_{t}=\tilde{\zeta}(t)\ \partial_{z}f_{t}
\end{equation}
has a canonical quasiconformal solution $f:\hat{G}\times \Delta \rightarrow \hat{G}\times \Delta$ such that $f(t,x):=\left(t,f_{t}(x)\right)$ if and only if $\zeta$ is continuous.
\end{lema}
\dem
If $f$ is quasiconformal then the map $t\mapsto f_{t}$ is continuous respect to the uniform convergence topology and because $\zeta(t)= \partial_{\bar{z}}f_{t}/\partial_{z}f_{t}$ then $\zeta$ is so. Conversely, by the Ahlfors-Bers Theorem, for every $t$ there is a unique canonical quasiconformal solution $f_{t}:\Delta \rightarrow \Delta$ to the equation \eqref{Ahlfors_Bers} such that the map $t\mapsto f_{t}$ is continuous respect to the uniform convergence topology.
\fdem

The laminated Poincar\'e disk is an unramified covering of the disk with boundary $\hat{G}\times S^{1}$. The case of a ramified laminated Poincar\'e disk with non trivial monodromy in its boundary is a much more difficult task. In \cite{BV}, we study the case of the laminated Poincar\'e disk with one point ramification with the rational adelic solenoid $\mathbb{A}_{\Q}/\Q$ as its boundary.

The previous Lemma puts constraints on the Beltrami differentials for the existence of deformations:

\begin{defi}
Consider the laminated uniformization $\pi:\hat{G}\times \Delta\rightarrow \Sigma_{\infty}$. We say that $\mu_{\infty}\in L_{\infty}^{cont}(\Sigma_{\infty},\eta\times m_{\infty})_{1}$ is a Beltrami differential on the lamination $\Sigma_{\infty}$ if $\pi^{*}(\mu_{\infty})$ is a continuous function $\hat{G}\rightarrow L_{\infty}(\Delta, m)_{1}$.
\end{defi}

By the previous Lemma, for every Beltrami differential $\mu$ on the lamination $\Sigma_{\infty}$ there is a unique canonical quasiconformal solution $f^{\mu}:\hat{G}\times \Delta\rightarrow \hat{G}\times \Delta$ to the corresponding Ahlfors-Bers equation:
$$\partial_{\bar{z}}f^{\mu}_{t}=\widetilde{\pi^{*}(\mu)}(t)\ \partial_{z}f^{\mu}_{t}$$
such that $f^{\mu}(t,x):=\left(t, f^{\mu}_{t}(x)\right)$. Define the $\mu$-deformed diagonal action as follows:
$$g\cdot_{\mu}(t,x):= \left(tg,\ f^{\mu}_{tg}\ \alpha(g)\ (f^{\mu}_{t})^{-1}\ (x)\right)$$

Define the $\mu$-deformed lamination $\Sigma_{\infty,\mu}$ as the quotient of the laminated Poincar\'e disk by the $\mu$-deformed diagonal action. Topologically, the deformed lamination is isotopic to the original one but with a different complex structure in general. Because of the equivariant relation:
$$f^{\mu}\left( g\cdot(t,x)\right)= g\cdot_{\mu} f^{\mu}(t,x)$$
we have a well defined quasiconformal map $[f^{\mu}]$ on $G$-orbits:
$$[f^{\mu}]:\Sigma_{\infty} \rightarrow \Sigma_{\infty,\mu}$$

Consider the set of pairs $(\Sigma_{\infty,\mu},[f^{\mu}])$ and define the following equivalence relation: We say that $(\Sigma_{\infty,\mu},[f^{\mu}])$ is \textit{Teichm\"uller equivalent} to $(\Sigma_{\infty,\nu},[f^{\nu}])$:
$$(\Sigma_{\infty,\mu},[f^{\mu}])\sim (\Sigma_{\infty,\nu},[f^{\nu}])$$
if there is a biholomorphism $h$ such that the following diagram commutes up to homotopy:

$$\xymatrix{ & \Sigma_{\infty} \ar[dl]_{[f^{\mu}]} \ar[dr]^{[f^{\nu}]} & \\
\Sigma_{\infty,\mu} \ar[rr]^{h} & & \Sigma_{\infty,\nu} }$$
Here a biholomorphism is a homeomorphism that leafwise it is a biholomorphism. Note that as a consequence of the definition, the biholomorphism $h$ must be leafwise preserving.

The space of equivalence classes $T(\Sigma_{\infty})$ has the structure of a fine moduli space. This is called the \textit{Teichm\"uller space} of the lamination $\Sigma_{\infty}$ and the universal family is called the \textit{Teichm\"uller universal curve}. The relation with the moduli $\mathcal{M}_{\infty}$ is the following\footnote{Considering the moduli spaces as stacks, the mapping class group is the fundamental group of the moduli space $\mathcal{M}$.}:
$$\mathcal{M}_{\infty}\cong T(\Sigma_{\infty})/MCG(\Sigma_{\infty})$$
where $MCG(\Sigma_{\infty})$ is the \textit{mapping class group} of the lamination:
$$MCG(\Sigma_{\infty}):=Homeo^{+}(\Sigma_{\infty})/Homeo_0(\Sigma_{\infty})$$

In \cite{Odden} it is proved that the mapping class group of the lamination is isomorphic to the group of virtual automorphisms of the fundamental group of the genus two surface $\Sigma$:
$$MCG(\Sigma_{\infty})\cong Vaut\left(\pi_{1}(\Sigma)\right)$$

\subsection{Limit-automorphic differentials}\label{Limit_automorphic_differentials}

\begin{defi}
Consider the representation of the group $G$ in the space of $L_{\infty}(\Delta)$ differentials isometries:
$$\rho:G\rightarrow Isom\left(L_{\infty}(\Delta)\right)$$
such that $\rho(g)=\alpha(g)^{*}$. A differential whose orbit is continuous respect to the profinite topology of $G$ will be called \textit{limit-automorphic}\footnote{To get familiar with the definition, suppose the group $G$ is the real line $\R$ acting in the space of functions $f:\R\rightarrow\R$ by translation: $(g\cdot f)(x)= f(x+g)$. Then, the functions whose orbits are continuous are precisely the uniform continuous ones.}. The space of limit-automorphic differentials will be denoted by $L_{\infty}(\Delta_{Emb})$.
\end{defi}

In concrete, a differential $\mu$ is limit-automorphic if for every $\varepsilon>0$ there is a finite index subgroup $G'$ of $G$ such that $||\alpha(g)^{*}(\mu)-\mu||_{\infty}<\varepsilon$ for every $g\in G'$. The space of limit-automorphic differentials $L_{\infty}(\Delta_{Emb})$ is a Banach subspace of the space of differentials $L_{\infty}(\Delta)$.

For every differential $\mu\in L_{\infty}(\Delta)$ define the function $\varphi_{\mu}:G\rightarrow L_{\infty}(\Delta)$ such that $\varphi_{\mu}(g):=\alpha(g)^{*}(\mu)$.

\begin{lema}\label{TFAE}
Consider a differential $\mu$ in $L_{\infty}(\Delta)$. Then, the following are equivalent:
\begin{enumerate}
\item $\mu$ is limit-automorphic.
\item $\varphi_{\mu}$ is continuous respect to the profinite topology of $G$. Moreover, its continuous extension $\hat{\varphi}_{\mu}(t)$ verifies:
\begin{equation}\label{G_equivariance_phi}
\varphi_{\mu}(gt)= \alpha(t)^{*}\varphi_{\mu}(g)
\end{equation}
for every $g$ in $G$ and every $t$ in $\hat{G}$.
\item There is a net of automorphic differentials $\left(\mu_{G'}\right)$ ranging over the finite index subgroups of $G$ such that $\mu_{G'}\in L_{\infty}(G')$ and the net uniformly converges to $\mu$.
\item There is a unique differential $\mu_{\infty}$ in $L_{\infty}^{cont}(\Sigma_{\infty},m_{\infty})$ such that:
$$(\pi\iota)^{*}(\mu_{\infty})= \mu$$
\end{enumerate}
\end{lema}
\dem
\textit{(1) implies (2):} The first assertion follows directly from the following identity:
$$||\varphi_{\mu}(g)-\varphi_{\mu}(h)||_{\infty}= ||\alpha(g)^{*}(\mu)-\alpha(h)^{*}(\mu)||_{\infty}=||\alpha(gh^{-1})^{*}(\mu)-\mu||_{\infty}$$
where we have used the fact that $||\gamma^{*}(\mu)||_{\infty}=||\mu||_{\infty}$ for every M\"obius transformation $\gamma$ and differential $\mu$.

For the second assertion we have:
$$\varphi_{\mu}(gt)= \alpha(gt)^{*}\mu= \alpha(t)^{*}\alpha(g)^{*}\mu= \alpha(t)^{*}\varphi_{\mu}(g)$$
for every pair of elements $g$ and $h$ in $G$. By continuity the result follows.

\textit{(2) implies (3):} Consider the unique continuous extension $\hat{\varphi}_{\mu}:\hat{G}\rightarrow L_{\infty}(\Delta)$ and define:
\begin{equation}\label{canonical_net}
\mu_{G'}:= [G':G]\ \int_{\overline{G'}}d\eta(t)\ \hat{\varphi}_{\mu}(t)
\end{equation}
where $G'$ is a finite index subgroup of $G$.

Every $\mu_{G'}$ is a $G'$-automorphic differential: Because of identity \eqref{G_equivariance_phi}, for every $g\in G'$ we have:
\begin{eqnarray*}
\alpha(g)^{*}(\mu_{G'}) &=& [G':G]\ \int_{\overline{G'}}d\eta(t)\ \alpha(g)^{*}\left(\hat{\varphi}_{\mu}(t)\right) = [G':G]\ \int_{\overline{G'}}d\eta(t)\ \hat{\varphi}_{\mu}(tg) \\
&=& [G':G]\ \int_{\overline{G'}g}d\eta(t)\ \hat{\varphi}_{\mu}(t) = [G':G]\ \int_{\overline{G'}}d\eta(t)\ \hat{\varphi}_{\mu}(t) = \mu_{G'}
\end{eqnarray*}

The net uniformly converges to $\mu$: By definition, for every $\varepsilon>0$ there is a finite index subgroup $S$ of $G$ such that $||\mu-\alpha(g)^{*}(\mu)||_{\infty}<\varepsilon$ for every $g\in S$. In particular, because of identity \eqref{G_equivariance_phi}, we have $||\mu-\hat{\mu}(g,\_)||_{\infty}<\varepsilon$ for every $g\in S$ and by continuity, we conclude that $||\mu-\hat{\varphi}_{\mu}(t)||_{\infty}<\varepsilon$ for every $t\in \overline{S}$. Then:
$$||\mu-\mu_{G'}||_{\infty}\leq [G':G]\ \int_{\overline{G'}}d\eta(t)\ ||\mu-\hat{\mu}(t,\_)||_{\infty}<\varepsilon\ [G':G]\ \int_{\overline{G'}}d\eta(t)= \varepsilon$$ 
for every finite index subgroup $G'$ of $S$ and we have the result.

\textit{(3) implies (1):} It follows directly from the following identity:
\begin{eqnarray*}
||\alpha(g)^{*}(\mu)-\mu||_{\infty} &=& ||\left(\alpha(g)^{*}(\mu)-\mu_{G'}\right)-\left(\mu-\mu_{G'}\right)||_{\infty} \\
&\leq& ||\alpha(g)^{*}\left(\mu-\mu_{G'}\right)||_{\infty}+||\mu-\mu_{G'}||_{\infty}= 2||\mu-\mu_{G'}||_{\infty}
\end{eqnarray*}
for every $g\in G'$ and every finite index subgroup $G'$ of $G$.

\textit{(2) implies (4):} Define the differential $\hat{\mu}$ in $L_{\infty}(\hat{G}\times \Delta,\eta\times m)$ such that $\hat{\mu}(t,z):=\hat{\varphi}_{\mu}(t)(z)$. Relation \eqref{G_equivariance_phi} implies that $\hat{\mu}$ factors through the quotient:
$$\Sigma_{\infty}= \left( \hat{G}\times \Delta\right)/G$$
by the diagonal action; i.e. there is a unique differential $\mu_{\infty}$ in $L_{\infty}(\Sigma_{\infty},m_{\infty})$ such that $\pi^{*}(\mu_{\infty})= \hat{\mu}$. Finally,
$$(\pi\iota)^{*}(\mu_{\infty})=\iota^{*}\pi^{*}(\mu_{\infty})= \iota^{*}\hat{\mu}= \mu$$

\textit{(4) implies (2):} Because of the diagonal action and the definition of $\mu_{\infty}$, the continuous function $\pi^{*}(\mu_{\infty})$ on $\hat{G}$ coincides with $\varphi_{\mu}$ on the dense subset $G$ and the result follows.
\fdem

The net defined in \eqref{canonical_net} will be called the \textit{canonical automorphic net} of the limit-automorphic differential $\mu$ and the differential $\hat{\mu}$ will be called the \textit{extended differential} of $\mu$.

\begin{lema}
Consider Beltrami differentials $\mu$ and $\nu$ on the lamination $\Sigma_{\infty}$. Then, $(\Sigma_{\infty,\mu},[f^{\mu}])\sim(\Sigma_{\infty,\nu},[f^{\nu}])$ if and only if the continuous extensions of $f^{\mu}$ and $f^{\nu}$ to the boundary $\hat{G}\times S^{1}$ coincide.
\end{lema}
\dem The proof is almost verbatim to Lemmas 5.1 and 5.2 of \cite{Imayoshi}.\fdem

As a corollary we have the following model for the Teichm\"uller space:
\begin{equation}\label{Model_Teichm}
T(\Sigma_{\infty})\cong L_{\infty}(\Delta_{Emb})_{1}/\sim
\end{equation}
where $\mu\sim\nu$ if the continuous extensions of $f^{\mu_{\infty}}$ and $f^{\nu_{\infty}}$ to the boundary $\hat{G}\times S^{1}$ coincide, with $\mu_{\infty}$ and $\nu_{\infty}$ as in item $(4)$ of the previous Lemma. This equivalence relation is the pullback of the equivalence relation upon Beltrami differentials in $L_{\infty}(\Delta)_{1}$ whose classes define the universal Teichm\"uller space $T(1)$. Because $L_{\infty}(\Delta_{Emb})$ is a Banach subspace of $L_{\infty}(\Delta)$, we have the closed embedding:
\begin{equation}\label{Teichm_embedding}
T(\Sigma_{\infty})\hookrightarrow T(1)
\end{equation}

\section{Renormalized Weil-Petersson metric}\label{Ren_Weil_Petersson}

The Weil-Petersson metric cannot be defined on the universal Teichm\"uller space\footnote{At the origin, it works only on tangent vectors whose corresponding vector fields on $S^{1}$ are $C^{3/2+\varepsilon}$. The corresponding vector field on $S^{1}$ of a differential $\mu$ is defined as follows:
$$\dot{f}[\mu]:=\frac{d\ f^{t\mu}|_{S^{1}}}{dt}\Bigr|_{t=0}$$}. Naively, one could hope that it could be defined on the Teichm\"uller space of the hyperbolic lamination but the situation gets worst: The metric evaluated on two limit-automorphic differentials is divergent or zero. A second approach could be to evaluate the Weil-Petersson metric on the canonical automorphic net of the respective vectors and take the limit. But this doesn't work either for the limit is divergent if the metric on the limit-automorphic differentials is so and we get stuck in the same place as before.

In this section we define a metric on the Teichm\"uller space of the lamination $\Sigma_{\infty}$ whose restricton to the Teichm\"uller space of the Riemann surface $\Sigma$ is the usual Weil-Petersson metric. We define the \textit{renormalized Weil-Petersson metric} on the tangent space at the origin $T_{0}T(\Sigma_{\infty})$ by the formula\footnote{We define it on representatives. Later we will see
 that this expression actually defines a metric on the tangent space of Teichm\"uller space; i.e. the space of infinitesimal Beltrami differentials equivalent to zero is the the space of Weil-Petersson null vectors.}:
\begin{equation}\label{Weil_Petersson}
(\mu, \nu)_{WP}:= \int_{\hat{G}}d\eta(t)\int_{F}dxdy\int_{\Delta}dx'dy'\ \frac{\hat{\mu}(t,z)\overline{\hat{\nu}(t,\zeta)}}{(1-z\bar{\zeta})^{4}}
\end{equation}
where $[\mu]$ and $[\nu]$ are limit-automorphic tangent vectors at the origin; i.e. $\mu, \nu\in L_{\infty}(\Delta_{Emb})$, and the differentials $\hat{\mu}$ and $\hat{\nu}$ are the $G$-invariant continuous extensions to $\hat{G}\times\Delta$ of the limit-automorphic differentials $\mu$ and $\nu$ respectively. Because of the $G$-invariance of the extended differentials and the $G$-equivariance of the trivial bundle $\hat{G}\times \Delta\rightarrow \Delta$, the definition \eqref{Weil_Petersson} is independent of the chosen fundamental domain. Because the Haar measure is invariant under the $G$-action, this metric is $G$-invariant:
\begin{equation}\label{invariance_Weil_Petersson}
(\alpha(g)^{*}\mu, \alpha(g)^{*}\nu)_{WP}=(\mu, \nu)_{WP}
\end{equation}
for every element $g$ in the group $G$.

The pullback of the complex structure on the universal Teichm\"uller space gives a complex structure on the Teichm\"uller space of the lamination $\Sigma_{\infty}$ and the renormalized Weil-Petersson metric \eqref{Weil_Petersson} is a K\"ahler metric respect to this complex structure.

As it was expected from the identity \eqref{Limit_Teichm} and the transversality result in \cite{NV}, the Teichm\"uller space of the lamination $\Sigma_{\infty}$ is transversal to the complex analytic and K\"ahler submanifold $Diff^{+}(S^{1})/PSL(2,\R)$. We reproduce their argument here adapted to our case.

Recall the Weil-Petersson metric $g_{WP}$ on the universal Teichm\"uller space:
$$g_{WP}(\mu, \nu) := \int_{\Delta}dxdy\int_{\Delta}dx'dy'\ \frac{\mu(z)\overline{\nu(\zeta)}}{(1-z\bar{\zeta})^{4}}$$

\begin{lema}\label{zero_or_infinity}
Consider limit-automorphic differentials $\mu$ and $\nu$ in $L_{\infty}(\Delta_{Emb})$. Then:
$$g_{WP}(\mu, \nu)=\sum_{g\in G}\ (\mu, \nu)_{WP}$$
In particular, $g_{WP}(\mu, \nu)$ is zero (diverges) if and only if $(\mu, \nu)_{WP}$ is zero (is not zero).
\end{lema}
\dem
The metric $g_{WP}$ is M\"obius invariant; i.e. For every M\"obius transformation $\gamma$ we have:
$$g_{WP}(\gamma^{*}\mu, \gamma^{*}\nu)=g_{WP}(\mu, \nu)$$
In particular, the metric is $G$-invariant and by Lemma \ref{TFAE} the map $\left( t\mapsto g_{WP}\left(\hat{\varphi_\mu}(t), \hat{\varphi_\nu}(t)\right)\right)$ is the constant $g_{WP}(\mu, \nu)$. Then:
\begin{eqnarray*}
g_{WP}(\mu, \nu) &=& \int_{\hat{G}}d\eta(t)\ g_{WP}(\mu, \nu) = \int_{\hat{G}}d\eta(t)\ g_{WP}\left(\hat{\varphi_\mu}(t), \hat{\varphi_\nu}(t)\right) \\
&=& \int_{\hat{G}}d\eta(t)\int_{\Delta}dxdy\int_{\Delta}dx'dy'\ \frac{\hat{\mu}(t,z)\overline{\hat{\nu}(t,\zeta)}}{(1-z\bar{\zeta})^{4}} \\
&=& \sum_{g\in G}\ \int_{\hat{G}}d\eta(t)\int_{g(F)}dxdy\int_{\Delta}dx'dy'\ \frac{\hat{\mu}(t,z)\overline{\hat{\nu}(t,\zeta)}}{(1-z\bar{\zeta})^{4}} \\
&=& \sum_{g\in G}\ \left(\alpha(g)^{*}(\mu), \alpha(g)^{*}(\nu)\right)_{WP} = \sum_{g\in G}\ (\mu, \nu)_{WP}
\end{eqnarray*}
Because the Fuchsian group is infinite, we have the result.
\fdem

\begin{rmrk}
\ \\
\begin{itemize}
\item As an immediate consequence of Lemma \ref{zero_or_infinity}, the null space respect to the Weil-Petersson metric $g_{WP}$ coincides with the one respect to the renormalized Weil-Petersson metric. In particular, the set of differentials equivalent to zero equals the renormalized Weil-Petersson null space:
$$\mu\sim 0\ \ \ iff\ \ \ \mu\in Null_{WP}\subset L_{\infty}(\Delta_{Emb})$$
Recall that, because this metric is positive semidefinite in $L_{\infty}(\Delta_{Emb})$, by Cauchy-Schwarz inequality the null space equals the set of null vectors.

\item The renormalized Weil-Petersson metric is not induced by the universal Weil-Petersson metric; i.e. the topological embedding $T(\Sigma_{\infty})\hookrightarrow T(1)$ is not isometric.
\end{itemize}
\end{rmrk}

\begin{teo}\label{Trans_Nag_Verj}
The Teichm\"uller space $T(\Sigma_{\infty})$ is transversal to the complex analytic submanifold $Diff^{+}(S^{1})/PSL(2,\R)$ in the universal Teichm\"uller space $T(1)$.
\end{teo}
\dem
Consider limit-automorphic differentials $\mu$ and $\nu$ in $L_{\infty}(\Delta_{Emb})$ such that $(\mu, \nu)_{WP}\neq 0$. By Lemma \ref{zero_or_infinity}, $g_{WP}(\mu, \nu)$ must diverge. In particular, if the corresponding vector field on $S^{1}$ of $\nu$ is $C^{\infty}$, then the one corresponding to $\mu$ cannot be $C^{3/2+\varepsilon}$ hence cannot be $C^{\infty}$. We have proved the following: If $\mu$ is a limit-automorphic differential whose corresponding vector field on $S^{1}$ is $C^{\infty}$, then $[\mu]=0$ for $||\mu||_{WP}=0$ and null vectors correspond to vectors equivalent to zero. The proof is complete.
\fdem

As it is explained in \cite{NV}, the previous result is an infinitesimal form of Mostow rigidity. Now, we explore the relation between the renormalized Weil-Petersson metric and the usual one. In what follows denote by $WP_{G'}$ the usual Weil-Petersson metric on $L_{\infty}(G')$.

\begin{lema}\label{relation_between_WPs}
Consider a finite index subgroup $G'$ of $G$ and $G'$-invariant differentials $\mu$ and $\nu$ in $L_{\infty}(G')$. Then:
$$(\mu, \nu)_{WP} = \frac{1}{[G':G]}\ WP_{G'}(\mu, \nu)$$
\end{lema}
\dem
Choose representatives $g_{i}\in G$ of the coset space $G'\backslash G$ such that the following union is a domain in the disk $\Delta$:
$$F':=\bigcup_{i=1}^{[G':G]}\ \alpha(g_{i})(F)$$
Notice that $F'$ is a fundamental domain of the subgroup $G'$. By definition \ref{Weil_Petersson} and linearity, we have the following decomposition:
\begin{eqnarray*}
(\mu, \nu)_{WP}  &:=& \int_{\hat{G}}d\eta(t)\int_{F}dxdy\int_{\Delta}dx'dy'\ \frac{\hat{\mu}(t,z)\overline{\hat{\nu}(t,\zeta)}}{(1-z\bar{\zeta})^{4}} \\
&=& \sum_{i=1}^{[G':G]}\int_{\overline{G'}g_{i}}d\eta(t)\int_{F}dxdy\int_{\Delta}dx'dy'\ \frac{\hat{\mu}(t,z)\overline{\hat{\nu}(t,\zeta)}}{(1-z\bar{\zeta})^{4}}
\end{eqnarray*}
By the $G$-equivariance \eqref{G_equivariance_phi} of the extensions we have:
\begin{eqnarray*}
(\mu, \nu)_{WP}  &=& \sum_{i=1}^{[G':G]}\int_{\overline{G'}}d\eta(t)\int_{F}dxdy\int_{\Delta}dx'dy'\ \frac{\alpha(g_{i})^{*}\hat{\mu}(t,z)\ \overline{\alpha(g_{i})^{*}\hat{\nu}(t,\zeta)}}{(1-z\bar{\zeta})^{4}} \\
&=& \sum_{i=1}^{[G':G]}\int_{\overline{G'}}d\eta(t)\int_{\alpha(g_{i})(F)}dxdy\int_{\Delta}dx'dy'\ \frac{\hat{\mu}(t,z)\overline{\hat{\nu}(t,\zeta)}}{(1-z\bar{\zeta})^{4}}
\end{eqnarray*}
Because $\mu$ and $\nu$ are $G'$-invariant, by \eqref{G_equivariance_phi} and continuity we have that the extension $\left(t \mapsto \hat{\mu}(t,z)\right)$ is constant on $\overline{G'}$ and equals $\mu(z)$. A similar result holds for the extension $\hat{\nu}$. Then:
\begin{eqnarray*}
(\mu, \nu)_{WP}  &=& \sum_{i=1}^{[G':G]}\int_{\overline{G'}}d\eta(t)\int_{\alpha(g_{i})(F)}dxdy\int_{\Delta}dx'dy'\ \frac{\mu(z)\overline{\nu(\zeta)}}{(1-z\bar{\zeta})^{4}} \\
&=& \int_{\overline{G'}}d\eta(t)\int_{F'}dxdy\int_{\Delta}dx'dy'\ \frac{\mu(z)\overline{\nu(\zeta)}}{(1-z\bar{\zeta})^{4}} \\
&=& \eta\left(\overline{G'}\right)\ WP_{G'}(\mu, \nu)
\end{eqnarray*}
Because of the $G$-invariance of the Haar measure, we have that $\eta\left(\overline{G'}\right)= [G':G]^{-1}$ and the proof is complete.
\fdem

As an immediate corollary, we have that the renormalized Weil-Petersson metric is an extension of the usual Weil-Petersson metric on $G$-invariant differentials:

\begin{cor}
Consider $G$-invariant differentials $\mu$ and $\nu$ in $L_{\infty}(G)$. Then:
$$(\mu, \nu)_{WP} = WP_{G}(\mu, \nu)$$
\end{cor}

\begin{cor}
Consider limit-automorphic differentials $\mu$ and $\nu$ in $L_{\infty}(\Delta_{Emb})$ such that $g_{WP}(\mu,\nu)\neq 0$; i.e. $g_{WP}(\mu,\nu)= \infty$ (Recall Lemma \ref{zero_or_infinity}). Then, the limit of the following net diverges:
$$\lim_{\substack{ \longleftarrow\\ G'<G \\ [G':G]<\infty}}\ WP_{G'}(\mu_{G'}, \nu_{G'}) = \infty$$
\end{cor}
\dem
By Lemma \ref{zero_or_infinity}, if $g_{WP}(\mu,\nu)\neq 0$ then $(\mu,\nu)_{WP}\neq 0$ and by definition it is finite. By continuity respect to the essential supremum norm and Lemma \ref{relation_between_WPs}, we have:
$$\lim_{\substack{ \longleftarrow\\ G'<G \\ [G':G]<\infty}}\ WP(\mu_{G'}, \nu_{G'})= \lim_{\substack{ \longleftarrow\\ G'<G \\ [G':G]<\infty}}\ [G':G]\ (\mu_{G'}, \nu_{G'}) = \infty$$
\fdem

By Lemma \ref{relation_between_WPs} and the continuity of the Weil-Petersson metric \eqref{Weil_Petersson}, we have the following characterization of the renormalized Weil-Petersson metric in terms of the usual one:

\begin{prop}
Consider the limit-automorphic differentials $\mu$ and $\nu$ in $L_{\infty}(\Delta_{Emb})$. Then:
$$(\mu, \nu)_{WP} = \lim_{\substack{ \longleftarrow\\ G'<G \\ [G':G]<\infty}}\ \frac{1}{[G':G]}\ WP_{G'}(\mu_{G'}, \nu_{G'})$$
\end{prop}

The previous Proposition justifies the name of the metric.

\begin{prop}\label{Limit_Teichm}
$$T\left(\Sigma_{\infty}\right)= \overline{\lim_{\substack{\longrightarrow \\ G'<G \\ [G':G]<\infty}} T\left(\Sigma_{G'}\right)}$$
\end{prop}
\dem
By Lemmas \ref{zero_or_infinity} and \ref{relation_between_WPs} we have the following commutative diagram:

$$\xymatrix{L_{\infty}(G) \ar@{^{(}->}[r] & \ldots L_{\infty}(G')  \ar@{^{(}->}[r] & \ldots L_{\infty}(\Delta_{Emb}) \ar@{^{(}->}[r] & L_{\infty}(\Delta)_{1} \\
Null(WP_{G}) \ar@{^{(}->}[r] \ar@{^{(}->}[u] & \ldots Null(WP_{G'}) \ar@{^{(}->}[r] \ar@{^{(}->}[u] & \ldots Null_{WP} \ar@{^{(}->}[r] \ar@{^{(}->}[u] & Null(g_{WP}) \ar@{^{(}->}[u] }$$
with $G'$ a finite index subgroup of $G$.

In particular, the net of inclusions of Beltrami spaces:
$$L_{\infty}(G)_{1}\subset \ldots L_{\infty}(G')_{1}\subset\ldots L_{\infty}(\Delta_{Emb})_{1}\subset L_{\infty}(\Delta)_{1}$$
induces the net of inclusions of Teichm\"uller spaces:
$$T(\Sigma_{G})\subset \ldots T(\Sigma_{G'})\subset \ldots T(\Sigma_{\infty})\subset T(1)$$
where $T(\Sigma_{G'})$ is the Teichm\"uller space of the Riemann surface $\Sigma_{G'}$ defined in \eqref{Sigma_G_prime}.

Because the canonical automorphic net  of a limit-automorphic differential converges uniformly in $L_{\infty}(\Delta_{Emb})$, the union of the finite Teichm\"uller spaces is dense in $T(\Sigma_{\infty})$ and we have the result.
\fdem

Because every Teichm\"uller space $T(\Sigma_{G'})$ is finite dimensional, in contrast with the usual universal 
Teichm\"uller space, the Teichm\"uller space of the lamination $\Sigma_{\infty}$ is a separable space.

\section{Theorem A}

To motivate this section, consider the finite index subgroups $G'$ and $G''$ of $G$ such that $G''<G'$ and $[G'':G']=n$. Then, we have a finite $n$-covering $\Sigma_{G''}\rightarrow \Sigma_{G'}$ where the Riemann surfaces are defined by \eqref{Sigma_G_prime}. Because the covering is non ramified, their Euler characteristic verify: $\chi(\Sigma_{G''})= n\ \chi(\Sigma_{G'})$. In particular:
$$g(\Sigma_{G''})= n \left( g(\Sigma_{G'})-1 \right)+1$$
where $g$ denotes the genus. Then, we have the following identities\footnote{Recall that the complex dimension of the Teichm\"uller space of a genus $g\geq 2$ compact Riemann surface is $3g-3$.}:
$$\dim_{\C} T(\Sigma_{G''})= n\ \dim_{\C} T(\Sigma_{G'})= \dim_{\C} T(\Sigma_{G''})^{n}= \dim_{\C} C\left( G''\backslash G', T(\Sigma_{G'})\right)$$

See that, because the left hand side does not depend on the group $G'$, it behaves as a ``free variable''. Besides their dimension, Is there any relation between the spaces $T(\Sigma_{G''})$ and $C\left( G''\backslash G', T(\Sigma_{G'})\right)$? Moreover, taking the profinite limit on $G''$, Is there any relation between the spaces $T(\Sigma_{\infty})$ and $C\left( \hat{G'}, T(\Sigma_{G'})\right)$?

To alleviate notation, we will simply denote the surface $\Sigma_{G'}$ by $\Sigma$ and the group $G'$ by $G$. The complex and K\"ahler structure on the space $C\left( \hat{G}, T(\Sigma)\right)$ with the compact-open topology\footnote{Because $\hat{G}$ is compact, this topology is the uniform convergence topology induced by the supremum metric respect to the Weil-Petersson distance.} are the those inherited as a function space from the target space $T(\Sigma)$ with the usual Weil-Petersson metric. The complex and K\"ahler structure on the space $T(\Sigma_{\infty})$ are those inherited from the complex analytic embedding \eqref{Teichm_embedding} and the renormalized Weil-Petersson metric \eqref{Weil_Petersson} defined in the previous section. The following Theorem shows that the relation we are looking for is the best one we could expect:

\begin{teo}\label{main1}
There is a complex analytic K\"ahler isometric homeomorphism between $C\left( \hat{G}, T(\Sigma)\right)$ and $T(\Sigma_{\infty})$:
$$\xymatrix{\hat{f}:C\left( \hat{G}, T(\Sigma)\right) \ar[r]^{\ \ \ \ \ \ \simeq}  & T(\Sigma_{\infty})}$$
\end{teo}

Consider a fundamental domain $F\subset\Delta$ of the $G$-action on the disk. Because of the $G$-equivariance of the trivial bundle $\hat{G}\times \Delta\rightarrow \Delta$, the set $\hat{G}\times F \subset \hat{G}\times \Delta$ is a fundamental domain of the $G$-action on the bundle.

Define the map $f:C\left( \hat{G}, L_{\infty}(G)\right) \rightarrow L_{\infty}(\Delta_{Emb})$ such that $f(\xi):=\mu$ where $\mu$ is the unique differential with the following property: The functions $\xi$ and $\hat{\varphi}_{\mu}$ coincide on the fundamental domain $\hat{G}\times F$:
$$\xi|_{\hat{G}\times F}= \hat{\varphi}_{\mu}|_{\hat{G}\times F}$$

We define the Weil-Petersson metric on the space $C\left(\hat{G}, L_{\infty}(G)\right)$ as follows:
$$\left(\mu, \nu\right)_{WP}:=\int_{\hat{G}}d\eta(t)\ \left(\mu(t), \nu(t)\right)_{WP}$$
where the metric on the right hand side is the usual Weil-Petersson metric on the space $L_{\infty}(G)$\footnote{Recall that the Weil-Petersson metric \eqref{Weil_Petersson} on $L_{\infty}(\Delta_{Emb})$ extends the usual one on the space $L_{\infty}(G)$.}.

\begin{lema}\label{Lemma_main_th}
Respect to the Weil-Petersson and renormalized Weil-Petersson metrics respectively, the map $f: C\left(\hat{G}, L_{\infty}(G)\right)\rightarrow L_{\infty}(\Delta_{Emb})$ is an isometric isomorphism\footnote{It is also an isometric isomorphism respect to the essential supreme norms: $||\mu||_{\infty}=||f(\mu)||_{\infty}$.}:
$$\left(\mu, \nu\right)_{WP}= \left(f(\mu), f(\nu)\right)_{WP}$$
\end{lema}
\dem
To see that $f$ is an isomorphism, it is enough to construct an inverse. Consider a limit-automorphic differential $\mu$ in $L_{\infty}(\Delta_{Emb})$ and define the function $\xi_{\mu}:\hat{G}\rightarrow L_{\infty}(G)$ such that:
$$\xi_{\mu}(t):=\sum_{g\in G} \alpha(g)^{*}\left(\hat{\varphi}_{\mu}(t)\ \chi_{F}\right)$$
This function $\xi$ is continuous for $\hat{\varphi}_{\mu}$ is continuous and the following relation:
\begin{eqnarray*}
||\xi_{\mu}(t)-\xi_{\mu}(h)||_{\infty} &=& ||\sum_{g\in G} \alpha(g)^{*}\left((\hat{\varphi}_{\mu}(t)-\hat{\varphi}_{\mu}(h))\ \chi_{F}\right)||_{\infty} \\
&=& ||(\hat{\varphi}_{\mu}(t)-\hat{\varphi}_{\mu}(h))\ \chi_{F}||_{\infty}\leq ||\hat{\varphi}_{\mu}(t)-\hat{\varphi}_{\mu}(h)||_{\infty}
\end{eqnarray*}
The map $(\mu\mapsto\xi_{\mu})$ is the linear inverse we were looking for.

Consider arbitrary continuous functions $\xi$ and $\zeta$ in $C\left( \hat{G}, L_{\infty}(G)\right)$. Because of their $G$-invariance, for every $t\in \hat{G}$ we have:
$$\xi(t)= \sum_{g\in G}\ \alpha(g)^{*}\left(\xi(t)\ \chi_{F}\right)$$
and an analogous identity for $\zeta(t)$. Then,
\begin{eqnarray*}
\left(\xi, \zeta\right)_{WP} &:=& \int_{\hat{G}}d\eta(t)\ \left(\sum_{g\in G}\ \alpha(g)^{*}\left(\xi(t)\ \chi_{F}\right), \sum_{h\in G}\ \alpha(h)^{*}\left(\zeta(t)\ \chi_{F}\right)\right)_{WP} \\
&=& \sum_{g,h\in G}\ \int_{\hat{G}}d\eta(t)\ \left(\alpha(g)^{*}\left(\xi(t)\ \chi_{F}\right), \alpha(h)^{*}\left(\zeta(t)\ \chi_{F}\right)\right)_{WP} \\
&=& \sum_{g\in G}\ \int_{\hat{G}}d\eta(t)\ \left(\alpha(g)^{*}\left(\xi(t)\ \chi_{F}\right), \alpha(g)^{*}\left(\zeta(t)\ \chi_{F}\right)\right)_{WP}
\end{eqnarray*}
The last step is because the $G$-action on the disk is free hence the integrand is proportional to a Kronecker delta:
\begin{equation}\label{delta_Kronecker}
\left(\alpha(g)^{*}\left(\xi(t)\ \chi_{F}\right), \alpha(h)^{*}\left(\zeta(t)\ \chi_{F}\right)\right)_{WP}\ =\ \delta_{g,h}\ \left(\alpha(g)^{*}\left(\xi(t)\ \chi_{F}\right), \alpha(g)^{*}\left(\zeta(t)\ \chi_{F}\right)\right)_{WP}
\end{equation}

Because the Haar measure is translation invariant, we can make the following variable change in the integral:
\begin{eqnarray*}
\left(\xi, \zeta\right)_{WP} &=& \sum_{g\in G}\ \int_{\hat{G}}d\eta(t)\ \left(\alpha(g)^{*}\left(\xi(tg^{-1})\ \chi_{F}\right), \alpha(g)^{*}\left(\zeta(tg^{-1})\ \chi_{F}\right)\right)_{WP} \\
&=& \int_{\hat{G}}d\eta(t)\ \left(\sum_{g\in G}\ \alpha(g)^{*}\left(\xi(tg^{-1})\ \chi_{F}\right), \sum_{h\in G}\ \alpha(h)^{*}\left(\zeta(tg^{-1})\ \chi_{F}\right)\right)_{WP} \\
&=& \int_{\hat{G}}d\eta(t)\ \left(\sum_{g\in G}\ g^{*}\left(\xi\ \chi_{\hat{G}\times F}\right)(t), \sum_{h\in G}\ h^{*}\left(\zeta\ \chi_{\hat{G}\times F}\right)(t)\right)_{WP} \\
&=& \int_{\hat{G}}d\eta(t)\ \left(\hat{\varphi}_{f(\xi)}(t), \hat{\varphi}_{f(\zeta)}(t)\right)_{WP}
\end{eqnarray*}
where we have used \eqref{delta_Kronecker}. Because of the renormalized Weil-Petersson $G$-invariance \eqref{invariance_Weil_Petersson}, we have\footnote{By the second item in Lemma \ref{TFAE} we have:
$$\hat{\varphi}_{f(\xi)}(t)= \alpha(g)^{*}\left(\hat{\varphi}_{f(\xi)}(tg^{-1})\right)$$
for every $t$ in $\hat{G}$ and $g$ in $G$. In particular, $\hat{\varphi}_{f(\xi)}(g)= \alpha(g)^{*}\left(\hat{\varphi}_{f(\xi)}(e)\right)$ for every $g$ in $G$.}:
\begin{eqnarray*}
\left(\hat{\varphi}_{f(\xi)}(g), \hat{\varphi}_{f(\zeta)}(g)\right)_{WP} &=& \left(\alpha(g)^{*}\left(\hat{\varphi}_{f(\xi)}(e)\right),  \alpha(g)^{*}\left(\hat{\varphi}_{f(\zeta)}(e)\right)\right)_{WP} \\
&=& \left(\hat{\varphi}_{f(\xi)}(e), \hat{\varphi}_{f(\zeta)}(e)\right)_{WP} = \left(f(\xi), f(\zeta)\right)_{WP}
\end{eqnarray*}
for every $g$ in the group $G$. By continuity, the above map uniquely extends on $\hat{G}$ to the constant map $\left(f(\xi), f(\zeta)\right)_{WP}$ and the result follows:
\begin{eqnarray*}
\left(\xi, \zeta\right)_{WP} &=& \int_{\hat{G}}d\eta(t)\ \left(\hat{\varphi}_{f(\xi)}(t), \hat{\varphi}_{f(\zeta)}(t)\right)_{WP} = \left(f(\xi), f(\zeta)\right)_{WP}
\end{eqnarray*}
\fdem

{\em Proof of Theorem \ref{main1}: \;} Recall that the set of (limit-automorphic) differentials equivalent to zero equals the set of (renormalized) Weil-Petersson null vectors. Because of Lemmas \ref{zero_or_infinity} and \ref{Lemma_main_th}, these null-vector sets are mapped one to one and the linear map $f$ descends to a map $\hat{f}$ on the respective Teichm\"uller spaces. \fdem

Now, we can think of the Teichm\"uller space $T(\Sigma_{\infty})$ as the space of continuous functions over a non-archimedean space taking values on the finite dimensional Teichm\"uller space $T(\Sigma)$.

We can also think of Theorem \ref{main1} as a canonical\footnote{Canonical up to the choice of a fundamental domain of the $G$ action on the disk.} assignment of K\"ahler coordinates on the Teichm\"uller space of the lamination.

Consider a finite index subgroup $G'$ of the group $G$ and its corresponding covering $\Sigma_{G'}$ defined in \eqref{Sigma_G_prime}. The homeomorphism $\hat{f}$ maps the subspace of locally constant functions that are constant on the closure of the $G'$-cosets to the finite dimensional subspace $T(\Sigma_{G'})$ of $T(\Sigma_{\infty})$; i.e. The following diagram commutes:

$$\xymatrix{ C\left( \hat{G}, T(\Sigma)\right) \ar[rrr]^{\ \ \ \ \ \ \hat{f}\ \ \simeq}  & & & T(\Sigma_{\infty}) \\
T(\Sigma)^{n}\simeq C\left( G'\backslash G, T(\Sigma)\right) \ar[rrr]^{\ \ \ \ \ \ \simeq} \ar@{^{(}->}[u]  & & & T(\Sigma_{G'}) \ar@{^{(}->}[u]
}$$
where $n$ is the index of $G'$ in $G$. Because every continuous function on a compact space is the uniform limit of locally constants functions, this gives another proof of Proposition \eqref{Limit_Teichm}.

\section{Theorem B}

Consider the symmetric group $Sym_{G'\backslash G}$ of the coset set $G'\backslash G$ and its natural right action on $C\left( G'\backslash G, T(\Sigma)\right)$:
$$C\left( G'\backslash G, T(\Sigma)\right) \times Sym_{G'\backslash G}\rightarrow C\left( G'\backslash G, T(\Sigma)\right)$$
such that $(f\cdot \varphi)(a)= f\left(\varphi(a)\right)$. Because the pushforward of the Haar measure on $G'\backslash G$ is just the normalized counting measure, by Lemma \ref{relation_between_WPs} the symmetric group acts by Weil-Petersson isometries on $C\left( G'\backslash G, T(\Sigma)\right)$. Consider the alternating subgroup $Alt_{G'\backslash G}$ of the isometries that preserve orientation. Then:
\begin{equation}\label{Isom_Sym}
Alt_{G'\backslash G}<Isom_{WP}^{+}\left(C\left( G'\backslash G, T(\Sigma)\right)\right)
\end{equation}

By the Masur-Wolf Theorem\footnote{This is the analog of Royden's Theorem for the Weil-Petersson metric. It states that the mapping class group is isomorphic to the positively oriented Weil-Petersson isometries:
$$MCG(\Sigma)\cong Isom^{+}_{WP}\left( T(\Sigma)\right)$$} \cite{Masur_Wolf} and Lemma \ref{relation_between_WPs}, we have:
\begin{equation}\label{MCG_Calculation}
MCG(\Sigma_{G'})= Isom_{WP}^{+}\left(T(\Sigma_{G'})\right)= Isom_{WP}^{+}\left(C\left( G'\backslash G, T(\Sigma)\right)\right)
\end{equation}
By the Masur-Wolf Theorem again, the group $C\left( G'\backslash G, MCG(\Sigma)\right)$ is a subgroup of the group of orientation preserving Weil-Petersson isometries on $C\left( G'\backslash G, T(\Sigma)\right)$:
\begin{equation}\label{Isom_MCG}
C\left( G'\backslash G, MCG(\Sigma)\right)< Isom_{WP}^{+}\left(C\left( G'\backslash G, T(\Sigma)\right)\right)
\end{equation}

By identities \eqref{MCG_Calculation}, \eqref{Isom_Sym} and \eqref{Isom_MCG} we have the following Mapping class subgroup:
$$\left\langle Alt_{G'\backslash G},\ C\left( G'\backslash G, MCG(\Sigma)\right)\right\rangle < MCG(\Sigma_{G'})$$
In particular, we have the following branched covering of the moduli space:
\begin{eqnarray*}
\mathcal{M}(\Sigma_{G'}) &=& T(\Sigma_{G'})/MCG(\Sigma_{G'})\twoheadleftarrow C\left( G'\backslash G, \mathcal{M}(\Sigma)\right)/Alt_{G'\backslash G}
\end{eqnarray*}

In a similar fashion as the symmetric product, we define the \textit{alternating product} of a space $X$ as follows:
$$Alt^{n}(X):=X^{n}/Alt_{n}$$
We have proved the following:

\begin{teo}
We have the following discrete fiber holomorphic K\"ahler isometric branched covering between moduli spaces factoring through the alternating product:
$$\mathcal{M}_{g}^{n}\twoheadrightarrow \mathcal{M}_{n(g-1)+1}$$
\end{teo}

\begin{rmrk}
Some remarks on the previous result:
\begin{itemize}
\item The covering is holomorphic K\"ahler isometric in the orbifold sense; i.e. outside the orbifold points and branch regions.
\item There is no explicit description, at least not an easy one that we were aware of, of the covering in terms of the isomorphism classes of Riemann surfaces. For example, the simplest case is:
$$\mathcal{M}_{2}\times\mathcal{M}_{2}\twoheadrightarrow \mathcal{M}_{3}$$
Is there any explicit way of describing a genus three Riemann surface by two genus two surfaces? Do not think about nodal curves, we are not in the Deligne-Mumford boundary.

\item This result is not in contradiction with the rigidity results in \cite{Aramayona_Souto}: Every holomorphic map $\mathcal{M}_{g, k}\rightarrow\mathcal{M}_{g', p}$ such that $g\geq 6$ and $g'\leq 2g-2$ is either constant or a forgetful map. In particular, under the previous hypothesis $p\leq k$.
\end{itemize}
\end{rmrk}

\section{Siegel functions and K\"ahler potential}\label{Siegel_Functions}

The Siegel disk $D(n)$ is defined as the space of $n\times n$ square complex symmetric matrices $Z$ such that $1-Z^{\dag}Z>0$. This space is a complex bounded domain of dimension $n(n +1)/2$ and it is identified with the homogeneous space $Sp(2n)/U(n)$\footnote{The group $Sp(2n)$ is the matrix subgroup of elements $g=
  \left( {\begin{array}{cc}
   a & b \\
   \bar{b} & \bar{a} \\
  \end{array} } \right)$ such that:
$$a^{t}\bar{b} -b^{\dag}a=0,\ \ \ \ \ \ a^{t} \bar{a}-b^{\dag}b=1$$
where $a$ and $b$ are $n\times n$ square complex matrices.}. In effect, consider the action of a $Sp(2n)$ element $g=
  \left( {\begin{array}{cc}
   a & b \\
   \bar{b} & \bar{a} \\
  \end{array} } \right)$ on a Siegel matrix $Z$:
$$g\cdot Z:= (aZ+b)(\bar{b}Z+\bar{a})^{-1}$$
This action is transitive and the isotropy group of the origin is the subgroup with $b=0$ isomorphic to $U(n)$.

The Siegel disk $D(n)$ is a K\"ahler space with K\"ahler potential:
$$K(Z)=-\tr \log(1-Z^{\dag}Z)$$

The Teichm\"uller space $T(\Sigma_{g})$ can be seen as the moduli space of marked Riemann surfaces of genus $g\geq 2$ with marking a set of generators of the fundamental group $\pi_{1}(\Sigma_{g})$. The abelianization map sends a set of generators of the fundamental group $\pi_{1}(\Sigma_{g})$ to a set of generators of the homology group $H_{1}(\Sigma_{g})$ hence a projection:
$$T(\Sigma_{g})\rightarrow TO(\Sigma_{g})$$
from the Teichm\"uller space to the respective Torelli space. Consider the canonical symplectic basis of cycles $\{a_{i},b_{j}\}$ such that:
$$a_{i} \cdot a_{j}= b_{i} \cdot b_{j} = \delta_{ij},\ \ \ \ \ a_{i}\cdot b_{j}=0$$
Consider the abelian differentials $\omega_{i}\in \Omega^{1}(\Sigma_{g})$ such that:
$$\delta_{ij} = \int_{a_{i}}\omega_{j}$$
and define the \textit{period matrix} $\Pi$ as follows:
$$\Pi_{ij}:= \int_{b_{i}}\omega_{j}$$
The period matrix $\Pi$ is a symmetric $n\times n$ complex matrix with positive imaginary part. Define the period mapping\footnote{The map $\xi$ is a holomorphic two to one immersion for $g\geq 3$ and a holomorphic embedding for $g=1,2$.} $\xymatrix{T(\Sigma_{g}) \ar[r]^{\pi} & TO(\Sigma_{g}) \ar[r]^{\xi} & D(g)}$ such that:
$$Z:= (\Pi-i)(\Pi+i)^{-1}$$
This is a holomorphic map and the pullback of the K\"ahler potential gives a K\"ahler potential on the Teichm\"uller space.

There is no direct analog of this construction for the universal Teichm\"uller space $T(1)$. However, there is a complex analytic embedding of $T(1)$ into the Segal disk $D_{Segal}$ constructed as follows: Consider the map $\varphi\to V_{\varphi}$ such that:
$$V_{\varphi}(f):= f\circ \varphi - \frac{1}{2\pi}\int_{S^{1}}f\circ \varphi$$
This map is introduced in \cite{NS} where it is proved that $V_{\varphi}$ is a bounded automorphism of the Hilbert space $\mathcal{H}:= H^{1/2}(S^{1},\C)/\C$ if and only if $\varphi$ is a quasisymmetric map of the circle\footnote{This map is a direct generalization of the Segal representation described in \cite{Segal}. It is also discussed in \cite{Nag}.}. The Hilbert space $\mathcal{H}$ has a symplectic structure given by the Hilbert transform and it decomposes in its $-i$ and $i$ eigenspaces $W_{+}$ and $W_{-}$ respectively: $\mathcal{H}= W_{+}\oplus W_{-}$. For every quasisymmetric map $\varphi$, the automorphism $V_{\varphi}$ preserves the symplectic structure on $\mathcal{H}$ and it is an isometry of the eigenspaces $W_{\pm}$ if and only if $\varphi$ is a M\"obius map. Hence we have the \textit{universal period mapping}:
\begin{equation}\label{UPM}
T(1):= QS(S^{1})/M\ddot{o}b\hookrightarrow D_{Segal}:= Sp(\mathcal{H})/U(W_{+})
\end{equation}

Composing with the embedding $T(\Sigma_{\infty})\hookrightarrow T(1)$, we have an embedding of the Teichm\"uller space $T(\Sigma_{\infty})$ in the Segal disk. This embedding is not satisfactory in the sense that the target space is too big with respect to the space $T(\Sigma_{\infty})$; i.e. The space $T(\Sigma_{\infty})$ is separable while the Segal disk is not. Looking forward to solve this problem, a first idea that comes to mind is to construct a direct limit of Siegel disks. However this strategy fails because there is no canonical way to embed the Siegel disk $D(n)$ into $D(n+1)$ and the lack of functoriality gives no map in the direct limit.

A solution to this problem is given by Lemma \ref{Lemma_main_th} for it gives the following holomorphic map:
\begin{equation}\label{holomorphic_map}
T(\Sigma_{\infty})\cong C\left( \hat{G}, T(\Sigma_{2})\right)\rightarrow C\left( \hat{G}, TO(\Sigma_{2})\right) \hookrightarrow C\left( \hat{G}, D(2)\right)
\end{equation}

The map \eqref{holomorphic_map} will be called the \textit{laminated period mapping} and elements of the function space $C\left( \hat{G}, D(2)\right)$ will be called \textit{Siegel functions}. The function space $C\left( \hat{G}, D(2)\right)$ has the following K\"ahler potential:
$$K(\zeta):= \int_{\hat{G}}\ d\eta(t)\ K\left(\zeta(t)\right)= -\int_{\hat{G}}\ d\eta(t)\ \tr \log\left(1-\zeta(t)^{\dag}\zeta(t)\right)$$
where $\zeta$ is a Siegel function. The pullback of this K\"ahler potential via the holomorphic map \eqref{holomorphic_map} gives a K\"ahler potential on the Teichm\"uller space $T(\Sigma_{\infty})$. The induced K\"ahler metric coincides with the one induced by the renormalized Weil-Petersson metric \eqref{Weil_Petersson}.

\section{Acknowledgments}
The first author benefited from C\'atedras CONACYT program.
The second author (AV) benefited from a PAPIIT (DGAPA, Universidad Nacional
Autónoma de México) grant IN106817.

\end{document}